# Vector bundles on a neighborhood of an exceptional curve and elementary transformations

## E. Ballico - E. Gasparim[*]

*Abstract. Let W be the germ of a smooth complex surface around an exceptional curve and let E be a rank 2 vector bundle on W. We study the cohomological properties of a finite sequence {$E_i$}$_{1 \leq i \leq t}$ of rank 2 vector bundles canonically associated to E. We calculate numerical invariants of E in terms of the splitting types of $E_i$, $1 \leq i \leq t$. If S is a compact complex smooth surface and E is a rank two bundle on the blow-up of S at a point, we show that all values of $c_2(E) - c_2(\pi_*(E)^{**})$ that are numerically possible are actually attained.*

## 0. Introduction

We consider exceptional curves in the following two cases. In the first case, let W be a smooth connected complex analytic surface which contains an exceptional divisor i.e. a smooth curve $D \cong \mathbf{P}^1$ with $\mathbf{O}_D(-1)$ as normal bundle. Let U be a small tubular neighborhood of D in the Euclidean topology and let $\pi: U \to Z$ be the contraction of D. In this case Z is the germ of a smooth surface around the point P$:= \pi$ (D).

In the second case, let W be a smooth connected algebraic surface defined over an algebraically closed field $\mathbf{K}$ with arbitrary characteristic. We assume that W contains an exceptional curve D and denote by U the formal completion of W along D. Let $\pi: U \to Z$ be the contraction of D. In this case Z is a formal smooth two-dimensional space supported at P.

In what follows we use the notation defined above to represent either case. Let $\mathbf{I}$ be the ideal sheaf of D in U and consider a rank 2 vector bundle E over U. Consider the pair of integers (a,b)  such that E|D $\cong \mathbf{O}_D(a) \oplus \mathbf{O}_D(b)$. We will refer to the pair (a,b) as the *splitting type* of E. Since Z is a smooth surface the bidual $\pi_*(E)^{**}$ is locally free and hence free because Z is 2-dimensional. There is a natural inclusion $\mathbf{j}: \pi_*(E) \to \pi_*(E)^{**}$ such that coker($\mathbf{j}$) has finite length. Set $Q:= $coker($\mathbf{j}$).

We show that the pair (z,w)$:= (h^0(Z,Q), h^0(Z,R^1\pi_*(E)))$ of numerical invariants of E is uniquely determined by a sequence of pairs of integers associated to E in [$\mathbf{B}$] using ele-mentary transformations. We review the construction of the associated sequence and prove the following results.

---
[*] The first author was partially supported by MURST (Italy) and the second author was partially supported by NSF grant DMS-0072675.

**Theorem 0.1.** *Let E be a rank 2 vector bundle on W with associated admissible sequence $\{(a_i, b_i)\}_{1 \le i \le t}$. Then we have the equalities*

$$w := h^0(Z, R^1 \pi_*(E)) = \Sigma_{1 \le i \le t} \max \{-b_i\text{-}1, 0\} \qquad and$$

$$z := h^0(Z, Q) = \Sigma_{1 \le i \le t} \, a_i - a_t^2 - \Sigma_{1 \le i \le t} \max \{-b_i\text{-}1, 0\}.$$

Every admissible sequence is associated to a rank 2 vector bundle on W (see [**B**] Th.0.2). For simplicity, we normalize our bundles to have splitting type $(j, -j+\varepsilon)$, with $\varepsilon = 0$ or $\varepsilon = -1$. We have the following existence theorem.

**Theorem 0.2.** *For every pair of integers $(z, w)$ satisfying $j\text{-}1\text{-}\varepsilon \le w \le j(j\text{-}1)/2 - j\varepsilon$ and $1 \le z \le j(j+1)/2$ with $j \ge 0$ and $\varepsilon = 0$ or $-1$, there exists a rank 2 vector bundle E on W with splitting type $(j, -j+\varepsilon)$ having numerical invariants $h^0(Z, R^1 \pi_*(E)) = w$ and $h^0(Z, Q) = z$.*

**Remark 0.3.** It follows from theorem 0.2 that the strata defined in [**BG**] for spaces of bundles on the blow-up of $C^2$ are all non-empty.

We give also the following characterization of the split bundle.

**Proposition 0.4.** *Let E be a rank 2 vector bundle on U with splitting type $(j, -j+\varepsilon)$ with $j > 0$ and $\varepsilon = 0$ or -1. The following conditions are equivalent:*

*(i)*    $E \cong O_U(-jD) \oplus O_U((j+\varepsilon)D)$

*(ii)*    $c_2(E) - c_2(\pi_*(E)^{**}) = j(j+\varepsilon)$

*(iii)*    $h^0(Z, R^1 \pi_*(E)) = j(j\text{-}1)/2$

*(iv)*    *E has associated sequence $\{(a_i, b_i)\}_{1 \le i \le j\text{-}\varepsilon}$ with $b_i = -j\text{-}\varepsilon+i\text{-}1$ for every i.*

We now consider a compact complex smooth surface S, so that we can calculate second chern classes. If E is a rank *2* bundle defined on the blow-up of S at a point, then the difference of second Chern classes satisfies $j \le c_2(E) - c_2(\pi_*(E)^{**}) \le j^2$ and is given by the sum $h^0(Z, R^1 \pi_*(E)) + h^0(Z, Q)$ (see [**FM**]). Sharpness of these bounds was proven in [**B**] and in [**G2**] by different methods. We prove the following existence theorem.

**Theorem 0.5.** *Let S be the blow-up of a compact complex smooth surface S at a point. Let l denote the exceptional divisor and let j be a non-negative integer. Then for every integer k satisfying $j \le k \le j^2$ there exists a rank 2 vector bundle* $E$ *over* $\mathbf{S}$ *with* $E|_l \cong \boldsymbol{O}_l(\mathrm{j}) \oplus \boldsymbol{O}_l(\text{-j})$ *satisfying* $c_2(E)$ - $c_2(\pi_*(E)^{**})$ =k.

**Note 0.6:** In [**G1**, Thm. 3.5] it is shown that the number of moduli for the space of rank-2 bundles on the blow up of $\mathbf{C}^2$ at the origin with splitting type j equals 2j-3; and since such bundles are determined by their restriction to a formal neighborhood of the exceptional divisor it follows that we have the same number of moduli for bundles over the neighborhood U of an exceptional curve on a surface W.

These results are proven in section 1, where we also review the construction of admissible sequences. On section 2 we consider briefly bundles of higher rank.

## 1. Rank 2 bundles

We briefly recall the construction of the associated sequences of pairs of bundles and splitting types given in the introduction of [**B**]. We first give the definitions of positive and negative elementary transformations.

Let E be a rank 2 vector bundle on W with splitting type (a,b) with $a \ge b$. Fix a line bundle R on D and a surjection **r**: E $\to$ R induced by a surjection $\rho$: E|D $\to$ R. There exists such a surjection if and only if deg(R) $\ge$ b. If deg(R) = b < a, then $\rho$ is unique, up to a multiplicative constant. Set $E' :=$ ker (**r**) and $R' =$ ker ($\rho$). If deg(R) = b < a the sheaf $E'$ is uniquely determined, up to isomorphism. Since D is a Cartier divisor, $E'$ is a vector bundle on U. We will say that $E'$ is the bundle obtained from E by making the *negative elementary transformation* induced by **r**. Note that $R'$ is a line bundle on D with degree deg($R'$) = a + b - deg(R). Since deg($\boldsymbol{I}/\boldsymbol{I}^2$) = 1 it is easy to check that deg($E'$|D) = a + b + 1 and we have the exact sequence

$$0 \to \boldsymbol{O}_D(1+\deg(R)) \to E'|D \to R' \to 0 \ . \tag{1}$$

Furthermore, using this exact sequence we obtain a surjection **t**: $E' \to R'$ such that ker(**t**) $\cong$ E(-D). In particular ker(**t**)|D $\cong \boldsymbol{O}_D(a+1) \oplus \boldsymbol{O}_D(b+1)$. Thus, up to twisting by $\boldsymbol{O}_U$(-D), the negative elementary transformation induced by **r** has an inverse operation and we will say that E is obtained from $E'$ making a *positive elementary transformation* supported by D. The following diagram, called the *display* of the elementary transformation, summarizes the construction (see [**M**]).

$$
\begin{array}{ccccccccc}
 & & 0 & & 0 & & & & \\
 & & \uparrow & & \uparrow & & & & \\
0 & \to & R' & \to & E|D & \xrightarrow{\ \rho\ } & R & \to & 0 \\
 & & {}_t\!\uparrow & & \uparrow & & \| & & \\
0 & \to & E' & \to & E & \xrightarrow{\ \mathbf{r}\ } & R & \to & 0 \\
 & & \uparrow & & \uparrow & & & & \\
 & & E(-D) & = & E(-D) & & & & \\
 & & \uparrow & & \uparrow & & & & \\
 & & 0 & & 0 & & & &
\end{array}
$$

Given two vector bundles $E_1$ and $E_2$ with splitting types $(a_1, b_1)$ and $(a_2, b_2)$ we say that $E_1$ is *more balanced* then $E_2$ if $a_1 - b_1 \le a_2 - b_2$. Given a vector bundles E with splitting type $(a, b)$ we say that E is *balanced* if either $a = b$ (case $c_1$ even) or else $a = b+1$ (case $c_1$ odd). Performing negative elementary transformations we will take the bundle E into more balanced bundles. The sequence of elementary transformations finishes when we arrive at a balanced bundle. If deg (R) $= b$, then $E'|D$ fits in the exact sequence

$$0 \to \boldsymbol{O}_D(b+1) \to E'|D \to \boldsymbol{O}_D(a) \to 0. \qquad (2)$$

If $b < a$ then $E'$ is more balanced than E. If $b \le a-3$, then (2) does not uniquely determine $E'|C$. If $b \le a-2$ and $E'$ is not balanced, we reiterate the construction starting from $E'$ taking $R'$ to be the factor of $E'|D$ of lowest degree and we take the unique surjection (up to a multiplicative constant) $\rho' : E' \to R'$. In a finite number, say, t-1, of steps, we send E into a bundle which, up to twisting by $\boldsymbol{O}_U(sD)$, where s= $-(a+b+t-1) / 2$ has trivial restriction to D. The process ends with a bundle isomorphic to $\boldsymbol{O}_U(sD)^{\oplus 2}$ (see [**B**], Remark 0.1).

We now construct the admissible sequence associated to E. Step one: set $E_1 := E$, $a_1 :=$ a and $b_1 :=$ b. If $a_1 = b_1$, set t = 1 and stop. Otherwise $a_1 > b_1$. Step two: in the case $a_1 > b_1$ set $E_2 := E'$ and let $(a_2,b_2)$ be the splitting type of $E'$. Note that $a_2+b_2 = a_1+b_1+1$ and $b_1 < b_2 \le a_2 \le a_1$. Hence $a_2 - b_2 < a_1 - b_1$ and $E_2$ is more balanced than $E_1$. If $a_2 = b_2$, set t:=2 and stop. If $a_2 > b_2$ reiterate the construction. Final step: in a finite number of steps (say t -1 steps) we arrive at a bundle $E_t$ with splitting type $(a_t,b_t)$ with $a_t = b_t$. Call $E_i$, $2 \le i \le t$, the bundle obtained after i-1 steps and let $(a_i,b_i)$ be the splitting type of $E_i$. The finite sequence of pairs $\{(a_i,b_i)\}_{1 \le i \le t}$ obtained in this way has the following properties:

    i.        $a_i \ge b_i \ \forall i > 0,$

    ii.      $a_i + b_i = a_1 + b_1 + i - 1 \ \forall i > 1,$

    iii.    $a_i \ge a_{i+1} \ge b_{i+1} > b_i \ \forall \ i \ge 1,$ and

    iv.    $a_t = b_t.$

We call *admissible* any such finite sequence of pairs of integers. We will say that a sequence $\{(a_i,b_i)\}, 1 \leq i \leq t$ is the *admissible sequence associated to the bundle E* if this sequence is created by the algorithm just described. By [**B**] Th. 0.2, every admissible sequence is associated to a rank 2 vector bundle on W.

*Examples:* Let us first set some notation. To represent the admissible sequence $\{(a_i,b_i)\}, 1 \leq i \leq t$, we write $(a_1,b_1) \to (a_2,b_2) \to \ldots \to (a_t,b_t)$.

1. If the splitting type of E is (b+2,b) then there is only one possibility for the admissible sequence associated to E, namely

$$(b+2,b) \to (b+2, b+1) \to (b+2, b+2).$$

2. If the splitting type of E is (b+4,b) then there are 3 different possibilities for admissible sequences associated to E (which in particular will give rise to different values of the numerical invariants (z,w) ), these are:

    i.   $(b+4,b) \to (b+4,b+1) \to (b+4,b+2) \to (b+4,b+3) \to (b+4, b+4)$
    ii.  $(b+4,b) \to (b+4,b+1) \to (b+3,b+3)$
    iii. $(b+4,b) \to (b+3,b+2) \to (b+3, b+3)$

We now calculate the numerical invariants of E in terms of admissible sequences. For every integer $n \geq 0$ let $D^{(n)}$ be the n-th infinitesimal neighborhood of D in U. Hence $D^{(n)}$ is the closed subscheme of U with $\mathbf{I}^{n+1}$ as ideal sheaf. In particular, $D^{(0)} = D$ and $D^{(n)}_{red} = D$ for every $n \geq 0$. For each integer $n \geq 0$ the following sequence is exact

$$0 \to \mathbf{I}^n/\mathbf{I}^{n+1} \to \boldsymbol{O}_U/\mathbf{I}^{n+1} \to \boldsymbol{O}_U/\mathbf{I}^n \to 0 \qquad (3)$$

Suppose that E is a vector bundle normalized to have splitting type $(j,-j+\varepsilon)$ where $j \geq 1$ and either $\varepsilon = 0$ or $\varepsilon = -1$. We denote by $\mathbf{m}$ be the maximal ideal of $\boldsymbol{O}_{Z,P}$. Consider the inclusion $\mathbf{j} \colon \pi_*(E) \to \pi_*(E)^{**}$ and let $Q := \operatorname{coker}(\mathbf{j})$, $z := h^0(Z,Q)$, and $w := h^0(Z,R^1\pi_*(E))$. Call $\boldsymbol{O}_D(x)$ the degree x line bundle on D. Twisting the exact sequence (3) by E and using the fact that $\mathbf{I}^n/\mathbf{I}^{n+1}$ has degree n, we obtain the exact sequence

$$0 \to \boldsymbol{O}_D(j+n) \oplus \boldsymbol{O}_D(-j+\varepsilon+n) \to E|D^{(n)} \to E|D^{(n-1)} \to 0 \qquad (4)$$

**Lemma 1.1.** *The integers z and w satisfy the inequalities:*

$$1 \leq z \leq j(j+1)/2 \quad and \quad j\text{-}1\text{-}\varepsilon \leq w \leq j\,(j\text{-}1)/2 - \varepsilon j.$$

*Proof.* By the Theorem on Formal Functions we have the bounds for z and we have that $w \leq \Sigma_{n \geq 0} \, h^1(D, \boldsymbol{O}_D(-j+\varepsilon+n)) = j(j-1)/2 - \varepsilon j$. The upper bound for w+z was stated in [**FM**] Remark 2.8, and proven for bundles with arbitrary rank in [**Bu**] Prop.2.8. Consequently we have an alternative proof of the upper bound for z. The lower bound for w will be proven in Remark 1.4. For the case of rank two and $\varepsilon = 0$ [**G2**] shows that these bounds are sharp.

Since $Q$ is a quotient of $\boldsymbol{O}_U, P^{\oplus 2}$ the dimension of the fiber of $Q$ at P is either 1 or 2. The sheaf $Q$ is isomorphic to the structure sheaf of a subscheme of Z supported by P and with length z if and only if the dimension of this fiber is 1. We will check that this is always true (see Proposition 1.3). We first check the split case.

**Lemma 1.2.** *Suppose that $E \cong \boldsymbol{O}_U(-jD) \oplus \boldsymbol{O}_U((j-\varepsilon)D)$ then $z = j(j+1)/2$, $w = j(j-1)/2 - \varepsilon j$ and Q is isomorphic to the structure sheaf of a subscheme of Z supported by P and with $\boldsymbol{m}^j$ as ideal sheaf.*

*Proof.* Since D is an exceptional divisor, we have $\pi_*(\boldsymbol{O}_U((j-\varepsilon)D)) = \pi_*(\boldsymbol{O}_U) = \boldsymbol{O}_Z$ for every $j \geq \varepsilon$ and $\pi_*(\boldsymbol{O}_U(-jD)) \cong \boldsymbol{m}^j$ if j>0.

**Proposition 1.3.** *Let E be a rank 2 vector bundle on W with splitting type $(j,-j+\varepsilon)$ with j > 0. Then Q is isomorphic to the structure sheaf of a length z subscheme Q of Z with $Q_{red} = P$ and $\boldsymbol{Q} \subseteq P^{(j-1)}$.*

*Proof of 1.3.* The first assertion is well-known and follows from the proof of Lemma 1.2. Since $Q$ is a quotient of $\boldsymbol{O}_Z{}^{\oplus 2}$, in order to prove the second assertion it is sufficient to check that its fiber at P is a 1-dimensional vector space. Since E has splitting type $(j,-j+\varepsilon)$, we have an extension

$$0 \to \boldsymbol{O}_U((-j+\varepsilon)\,D) \to E \to \boldsymbol{O}_U(\,jD) \to 0 \qquad (5)$$

([**BG**] Lemma 1.2, or in [**G1**] Thm. 2.1 in the case $\varepsilon = 0$ ). Call **e** the extension (5) giving E. For each $t \in \boldsymbol{K} \backslash \{0\}$ consider the extension of $\boldsymbol{O}_U(jD)$ by $\boldsymbol{O}_U((-j+\varepsilon)D)$ given by extension class t**e**, this extension has as middle term a vector bundle isomorphic to E. Using the extension **e** for t = 0, we construct a family $\{\lambda\boldsymbol{e}\}_{\lambda \in \boldsymbol{K}}$ of extensions. We call $E_\lambda$ the corresponding middle term and $Q_\lambda$ the corresponding sheaf. Since $E_\lambda \cong E$ for $\lambda \neq 0$, we have $Q_\lambda = Q$ for $\lambda \neq 0$, and because $E_0 \cong \boldsymbol{O}_W(jD) \oplus \boldsymbol{O}_W((-j+\varepsilon)D)$, we have that $Q_0 = P^{(j-1)}$, and the result follows by semi-continuity of the fiber dimension at P.

*Proof of 0.1.* Given the admissible sequence of splitting types $\{(a_i,b_i)\}_{1 \leq i \leq t}$ associated to E we want to show that

$$w := h^0(Z, R^1\pi_*(E)) = \Sigma_{1 \leq i \leq t} \max \{-b_i - 1, 0\} \qquad \text{and}$$
$$z := h^0(Z, Q) = \Sigma_{1 \leq i \leq t} \; a_i - a_t^2 - \Sigma_{1 \leq i \leq t} \max \{-b_i - 1, 0\}.$$

We use induction on t, the case t = 1 arising if and only if $a_1 = b_1$, equivalently, when $E \cong \boldsymbol{O}_W(-a_1 D)^{\oplus 2}$ (this follows immediately from the definition of admissible sequence). Since

$$R^1\pi_*(\boldsymbol{O}_W(xD)) = 0 \; \forall \; x \leq 1 \quad \text{and}$$
$$R^1\pi_*(\boldsymbol{O}_W(yD)) = y(y-1)/2 \; \; \forall \; y > 0,$$

we have the equality for w in the split case. Assume $t \geq 2$. By the definition of the sequence $\{E_i\}_{1 \leq i \leq t}$ associated to E we have that $E_1 = E$ and that there is an exact sequence

$$0 \to E_2 \to E_1 \to \boldsymbol{O}_D(-b_1 D) \to 0. \tag{6}$$

First assume $b_1 < 0$, in which case we have that $h^0(Z, \pi_*(\boldsymbol{O}_D(-b_1 D))) = 0$ and $h^0(Z, R^1\pi_*(\boldsymbol{O}_W(-b_1 D))) = -b_1 - 1$. Hence $w := h^0(Z, R^1\pi_*(E)) = h^0(Z, R^1\pi_*(E_2)) - b_1 + 1$ and since $E_2$ has $\{(a_{i+1}, b_{i+1})\}_{1 \leq i < t}$ as admissible sequence, the claim follows.

Now assume $b_1 \geq 0$, from the exact sequence (6) it follows that $h^0(Z, R^1\pi_*(E)) \leq h^0(Z, R^1\pi_*(E_2))$. Since $b_i > b_1$ for every $i > 1$, we have $h^0(Z, R^1\pi_*(E_2)) = 0$. Hence, by the inductive assumption on the length of the admissible sequence, it follows that $h^0(Z, R^1\pi_*(E)) = 0$, proving the first assertion. The value of $z := h^0(Z, Q)$ comes from the equalities $c_2(E) - c_2(\pi_*(E)^{**}) = \Sigma_{1 \leq i < t} \; a_i - a_t^2$ and $c_2(E) - c_2(\pi_*(E)^{**}) = h^0(Z, Q) + h^0(Z, R^1 \pi_*(E))$ proved in [**B**, Th. 0.3] and in [**FM**] respectively. Here, of course, we assume that E is extended to a compactification, however these integers do not depend upon the choices of compactification and of extension of E.

*Proof of 0.2.* By [**B**] Th. 0.2 every admissible sequence $(a_i, b_i)$ is associated to a rank two bundle E on W, moreover, the intermediate steps of the construction of E give bundles $E_i$ with splitting types $(a_i, b_i)$ for each i. Now use Th. 0.1 to calculate z and w.

**Remark 1.4.** If we assume that E has splitting type $(j, -j+\varepsilon)$ with $j \geq 1+\varepsilon$, then because $b_1 = -j+\varepsilon$, we obtain $w \geq j-1-\varepsilon$.

*Proof of 0.4.* By [**B**] Th. 0.5 we know that (i) and (ii) are equivalent. By Lemma 1.2 (i) implies (iv). Since $b_1 = -j + \varepsilon$, $b_i > b_{i-1} \; \forall \; i > 1$ holds, and $a_1 = j$ and $a_i + b_i = \varepsilon + i - 1 \; \forall \; i$, it follows from Theorem 0.1 that (iv) implies (ii).

*Proof of 0.5.* Given bundles G on S and F on W with $c_1(G)=0=c_1(F)$ there exists a bundle E on **S** satisfying E | **S** – $l$ =$\pi$* E| S-{p} and E|$W = F$ (see **[G3]** Cor. 3.4).It then follows that $c_2(E) - c_2(\pi_*(E)^{**}) = R^l\pi_*(F) + l(Q)$ and the result follows from Th. 0.2.

## 2. Bundles of higher rank

In this section we consider vector bundles with rank r ≥ 3. Fix a rank r vector bundle E on U. We use the notation of **[BG]** ß3 for the admissible sequence $\{E_i\}$, $_{1 \le i \le t}$ of vector bundles associated to E. In particular we denote by $(a(i,1),...,a(i,r))$ the splitting type of $E_i$ where for $a(i,1) \ge ... \ge a(i,r)$. We make the strong assumption that $a(i,r-1) \ge -1$ for every i and compute $h^0(Z,R^l\pi_*(E))$.

**Proposition 2.1.** *Let E be a rank r vector bundle on W whose associated sequence of vector bundles $\{E_i\}$ has splitting type $(a(i,1),..., a(i,r))$ with $1 \le i \le t$ and $a(i,r-1) \ge -1$, for all i. Then we have $h^0(Z,R^l\pi_*(E)) = \Sigma_{1\le i \le t} min\{-a(r,i)-1,0\}$.*

*Proof.* We first observe that the proof of the corresponding inequality for rank 2 bundles works verbatim (both cases t = 1 and t > 1), because for each i with $1 \le i \le t$ at most one of the integers a(i,j) is not at least –1 and $h^1(\mathbf{P}^1,L) = 0$ for every line bundle L on $\mathbf{P}^1$ with deg(L) ≥-1.

In the case r ≥ 3, the sequence of elementary transformations made to balance the bundle is not, a priori, uniquely determined, and hence the sequence of associated bundles is not uniquely determined by E. The condition a(1,r-1) ≥-1 implies that there is an associated sequence in which we make always an elementary transformation with respect to $\mathbf{O}_D(a(r,i))$ to pass from $E_i$ to $E_{i+1}$ for some a(r,i) ≤ -1 (which gives that $h^0(Z,R^l\pi_*(E_i))$ = $h^0(Z,R^l\pi_*(E_{i+1}))$- a( r, i ) + 1). We continue to perform elementary transformations until we arrive at an integer m ≤ t such that a(m,j) ≥ -1 for every i. It is then quite easy to check that $h^0(Z,R^l\pi_*(E_m)) = 0$ and the result follows.

In the general case the same method gives the following partial result.

**Proposition 2.2.** *Let E be a rank r vector bundle on W whose associated sequence $\{E_i\}$, $1 \le i \le t$ of vector bundles has splitting type $(a(i,1),...,a(i,r))$ with $1 \le i \le t$. Then we have $h^0(Z,R^l\pi_*(E)) \le \Sigma_{1 \le i \le t,\ 1 \le j \le r} min\{-a(j,i)-1,0\}$.*

Edoardo Ballico

Dept. of Mathematics, University of  Trento

38050 Povo (TN) - Italy

fa0:  39-*0461881624*

e-mail: ballico@science.unitn.it

Elizabeth Gasparim

Department of Mathematics

University of Texas at Austin

Austin TX 78712

e-mail:  gasparim@math.utexas.edu